\documentclass[11pt,english]{article}
\usepackage[LGR,T1]{fontenc}
\usepackage[latin9]{inputenc}
\usepackage{amsmath}
\usepackage{amsthm}
\usepackage{amssymb}

\makeatletter

\DeclareRobustCommand{\greektext}{%
  \fontencoding{LGR}\selectfont\def\encodingdefault{LGR}}
\DeclareRobustCommand{\textgreek}[1]{\leavevmode{\greektext #1}}
\ProvideTextCommand{\~}{LGR}[1]{\char126#1}

\theoremstyle{plain}
\newtheorem{thm}{\protect\theoremname}
  \theoremstyle{definition}
  \newtheorem*{problem*}{\protect\problemname}

\usepackage{mathrsfs}\usepackage{amsthm}\usepackage{amscd}

\usepackage{hyperref}
\usepackage{graphics}
\usepackage{a4wide}
\@ifundefined{definecolor}
 {\usepackage{color}}{}

\newcounter{EQNR}
\renewcommand{\contentsname}{Contents\\
{\footnotesize\normalfont(A table of contents should normally not
be included)}}

\makeatother

\usepackage{babel}
  \providecommand{\problemname}{Problem}
\providecommand{\theoremname}{Theorem}

\begin{document}

\title{Spectral zeta functions}

\author{Anders Karlsson \thanks{Preliminary version. October, 2018. Supported in part by the Swiss NSF.}}
\maketitle
\begin{abstract}
This paper discusses the simplest examples of spectral zeta functions,
especially those associated with graphs, a subject which has not been
much studied. The analogy and the similar structure of these functions,
such as their parallel definition in terms of the heat kernel and
their functional equations, are emphasized. Another theme is to point
out various contexts in which these non-classical zeta functions appear.
This includes Eisenstein series, the Langlands program, Verlinde formulas,
Riemann hypotheses, Catalan numbers, Dedekind sums, and hypergeometric
functions. Several open-ended problems are suggested with the hope
of stimulating further research.
\end{abstract}

\section{Introduction}

Euler observed the following product formula:
\[
\sum_{n=1}^{\infty}\frac{1}{n^{s}}=\prod_{p}(1-p^{-s})^{-1}
\]
where the product is taken over the prime numbers. This function of
the complex variable $s$ is called the Riemann zeta function, denoted
$\zeta(s)$, and the expressions above are convergent for $\mathrm{Re}(s)>1$.
The right hand side inspired several generalizations, by Artin, Hasse,
Weil, Selberg, Ihara, Artin, Mazur, Ruelle, and others, see \cite{IK04,Te10}
and references therein. The most far-reaching frameworks for Euler
products might be provided by the insights of Grothendieck and Langlands.
The left hand side found generalizations by Dirichlet, Dedekind, Hurwitz,
Epstein and others in number theory, but also in another direction
in the work of Carleman \cite{Ca34}, later extended in \cite{MP49},
where instead of the integers one takes Laplacian eigenvalues:
\[
\sum_{\lambda\neq0}\frac{1}{\lambda^{s}}
\]
convergent in some right half-plane. The purpose of the present note
is to survey a few recent investigations of these latter functions,
\emph{spectral zeta functions}, in cases where the Laplacian is less
classical, instead coming from graphs or p-adics. As in one of Riemann's
arguments, to get the analytic continuation one rather defines the
zeta function via the heat kernel, so the Mellin transform of the
heat trace (removing the constant term whenever necessary) and divide
by $\Gamma(s)$. This has the advantage of also making sense when
the spectrum has continuous part. 

One intriguing aspect that we mention is the functional equation of
the type $s$ vs $1-s$ that appears also for these non-classical
zeta functions. Notably one has
\[
\xi_{\mathbb{Z}}(1-s)=\xi_{\mathbb{Z}}(s),
\]
see below, and the equivalence of certain asymptotic functional equations
to the Riemann hypotheses for $\zeta(s)$ and certain Dirichlet $L$-functions
\cite{FK17,F16}.

Although there are a few instances in the literature where such function
are introduced for graphs, it seems that the first more systematic
effort to study spectral zeta functions of graphs appear in my paper
with Friedli \cite{FK17}. As it turns out, via asymptotic considerations,
these functions are intimately related to certain zeta functions from
number theory. In what follows we will moreover point out that spectral
zeta functions for graphs appear \emph{incognito} in, or are connected
to, a rich set of topics:
\begin{itemize}
\item Eisenstein series, continuous spectrum of surfaces, Langlands program, 
\item Riemann hypotheses 
\item Dedekind sums and Verlinde formulas, 
\item hypergeometric functions of Appell and Lauricella
\item Catalan numbers
\item Fuglede-Kadison determinants
\end{itemize}
Although each of these connections may not be of central importance
for the corresponding topic, still, in view of the variety of such
appearances, I think that the subject of graph spectral zeta functions
deserves more attention. This note can also be viewed as a modest
update to parts of the discussion in Jorgenson-Lang \cite{JL01} originally
entitled \emph{Heat kernels all over the place}. Throughout the text,
I will suggest some questions and further directions for research.

I would like the thank the organizers, Mugnolo, Atay and Kurasov,
for the very stimulating workshop \emph{Discrete and continuous models
in the theory of networks, }at the ZIF, Bielefeld, 2017, and for their
invitation. I thank Fabien Friedli, Mårten Nilsson and the referee
for corrections.

\section{Continuous case}

\subsection{$\mathbb{R}$ and $\mathbb{R\mathrm{/\mathbb{Z}}}$ }

This case has been recorded in so many places and there is no need
to repeat it here. But for comparison with the other contexts, we
recall the formulas. The heat kernel on $\mathbb{R}$ is 
\[
K_{\mathbb{R}}(t,x)=\frac{1}{\sqrt{4\pi t}}e^{-x^{2}/t}.
\]
 The heat kernel on the circle $\mathbb{R\mathrm{/\mathbb{Z}}}$ is
\[
K_{\mathbb{R\mathrm{/\mathbb{Z}}}}(t,x)=\sum_{n\in\mathbb{Z}}e^{-4\pi^{2}n^{2}t}e^{2\pi inx}.
\]
 This expression comes from spectral considerations and equals the
following periodization (i.e. one sums over the discrete group $\mathbb{Z}$
to obtain a function on the quotient $\mathbb{R\mathrm{/\mathbb{Z}}}$,
like in the proof of the Poisson summation formula) :
\[
\frac{1}{\sqrt{4\pi t}}\sum_{m\in\mathbb{Z}}e^{-(x+n)^{2}/4t}.
\]
 Setting $x=0$ in the former expression, removing $1,$ and taking
the Mellin transform and dividing by $\Gamma(s)$ yields
\[
\zeta_{\mathbb{R\mathrm{/\mathbb{Z}}}}(s):=2\cdot4^{-s}\pi^{-2s}\zeta(2s),
\]
which is the spectral zeta function of the circle $\mathbb{R\mathrm{/\mathbb{Z}}}$.
If we replace $s$ by $s/2$ we have
\[
2\cdot(2\pi)^{-s}\zeta(s).
\]
By defining the completed zeta function (in particular by multiplying
back the gamma factor)
\[
\xi(s):=\frac{1}{2}2^{s}\pi^{s/2}\Gamma(s/2)\zeta_{\mathbb{R\mathrm{/\mathbb{Z}}}}(s/2)
\]
and equating the two heat kernel expressions (Poisson summation formula)
one gets the fundamental functional equation in this form 
\[
\xi(1-s)=\xi(s),
\]
known in the physics literature as the reflection formula for $\zeta(s)$. 

The identity with the two heat kernel expressions are known to contain
moreover a wealth of theorems, such as the modularity of theta functions
and the law of quadratic reciprocity. 

\subsection{A few comments on further examples}

As recalled in the introduction the definition of spectral zeta function
was perhaps first in \cite{Ca34}, in fact he had more general functions
using also the eigenfunctions. This was used for the study of asymptotic
properties of the spectrum using techniques from analytic number theory.
Another application is for the purpose of defining determinants of
Laplacians, in topology this was done by Seeley, Ray and Singer in
the definition of analytic torsion, see \cite{R97}, and in physics
by Dowker, Critchley, and Hawking \cite{DC76,Ha77}, useful in several
contexts, see \cite{E12}. A referee suggested that in this context
one can also make reference to the important heat kernel approach
to index theorems, see \cite{BGV92}, and to the connection with Arakelov
geometry, see \cite{So92}.

The determinant of the Laplacian is defined in the following way,
assuming the analytic continuation of the zeta function to $s=0$,
\[
\det\Delta_{X}:=e^{-\zeta_{X}'(0)}.
\]

For the numbers $1,2,3,...$ (essentially the circle spectrum) the
corresponding determinant, which formally would be $1\cdot2\cdot3...=\infty!$,
has the value $\sqrt{2\pi}$ thanks to the corresponding well-known
special value of $\zeta(s)$, and this value is coherent with the
asymptotics of the factorial function $n!$ by de Moivre and Stirling.
In fact it was exactly this constant that Stirling determined. 

For tori the corresponding spectral zeta function are Epstein zeta
functions. There are also some studies of spheres with explicit formulas.
Osgood, Phillips, Sarnak \cite{OPS88} showed that the determinant
of the Laplacian is a proper function on the moduli spaces of Riemann
surfaces, and could conclude that the set of isospectral surfaces
is pre-compact. In string theory it was important to study how the
determinant of laplacian changes when varying the metric on surfaces. 

Spectral functions of Riemannian manifolds will appear in the limit
of spectral zeta functions for certain sequences of graphs. The first
more substantial example, the case of tori, can be found in \cite{FK17}.

\section{Discrete case}

To avoid confusion, let me right away emphasize that I am not considering
the Ihara zeta function or any related function. Therefore I do not
give references for this zeta function, other than the book by Terras
\cite{Te10} and another recent references for the Ihara zeta function
of quantum graphs \cite{Sm07}. 

On the other hand, what is more relevant in our context here, but
will not be discussed, are the spectral zeta functions of quantum
metric graphs studied by Harrison, Kirsten and Weyand, see \cite{HK11,HWK16,HW18}
and references therein. The paper \cite{HW18} even discusses the
spectral zeta functions of discrete graphs exactly in our sense.

Other related papers on spectral functions of (especially finite cyclic)
graphs are those of Knill \cite{Kn13,Kn18} which contain a wealth
of interesting ideas, and some considerations close to topics in \cite{FK17}
(we do not however understand his Theorem 1c in \cite{Kn13}, the
convergence to a non-vanishing function in the critical strip, see
also his Theorem 10; for us this limiting function is the Riemann
zeta function which of course has zeros.) Knill also defines a function
$c(s)$ that coincides with our $\zeta_{\mathbb{Z}}(s)$ perhaps without
connecting it to the graph $\mathbb{Z}$. The asymptotics for cyclic
graphs considered in \cite{Kn13} and \cite{FK17} were also deduced
in \cite{Si04}. The re-proof of Euler's formulas for $\zeta(2n)$
via such asymptotics also appears in \cite{CJK10}, but also here
there are earlier references, for example \cite{W91}, in view of
the form of Verlinde's formulas, discussed below. 

\subsection{$\mathbb{Z\mathrm{}}$ and $\mathbb{Z\mathrm{/}}n\mathbb{Z}$}

Using the heat kernel on $\mathbb{Z}$ in terms of the $I$-Bessel
function, see for example \cite{KN06} where it is rediscovered, the
spectral zeta function of the graph $\mathbb{Z}$ is
\begin{equation}
\zeta_{\mathbb{Z}}(s)=\frac{1}{\Gamma(s)}\int_{0}^{\infty}e^{-2t}I_{0}(2t)t^{s}\frac{dt}{t},\label{eq:zetaZ}
\end{equation}
where it converges as it does for $0<Re(s)<1/2$. From this definition
it is not immediate that it admits a meromorphic continuation and
a functional equation very similar to classical zetas. However, the
following was shown in \cite{FK17}:
\begin{thm}
For all $s\in\mathbb{C}$ it holds that 
\[
\zeta_{\mathbb{Z}}(s)=\frac{1}{4^{s}\sqrt{\pi}}\frac{\Gamma(1/2-s)}{\Gamma(1-s)}.
\]
The function
\[
\xi_{\mathbb{Z}}(s)=2^{s}\cos(\pi s/2)\zeta_{\mathbb{Z}}(s/2),
\]
is entire and satisfies for all $s\in\mathbb{C}$ 
\[
\xi_{\mathbb{Z}}(s)=\xi_{\mathbb{Z}}(1-s).
\]
\end{thm}

Jérémy Dubout \cite{Du16} observed that in fact one can write this
function as follows:
\begin{equation}
\zeta_{\mathbb{Z}}(s)=\left(\begin{array}{c}
-2s\\
-s
\end{array}\right).\label{eq:Dubout-formula}
\end{equation}
This makes a connection to the \textbf{Catalan numbers} 
\[
C_{n}=\frac{1}{n+1}\left(\begin{array}{c}
2n\\
n
\end{array}\right)
\]
which are ubiquitous in combinatorics (214 such manifestations are
listed in the book by Stanley \cite{St12}). 
\begin{problem*}
Could the other $\zeta_{\mathbb{Z}^{d}}(-n)$ be thought of as generalizations
of the Catalan numbers?
\end{problem*}
Note also that (\ref{eq:Dubout-formula}) immediately shows, what
is not a priori clear from the definition (\ref{eq:zetaZ}), that
at negative integers this zeta function takes rational (indeed integral)
values. This is analogous to Riemann zeta function and other Dedekind
zeta functions, by theorems of Hecke, Klingen, Siegel, Deligne and
Ribet. I refer to \cite{Du16} for a fuller discussion on this topic
in the graph setting. 

\textbf{Eisenstein series} are functions which appear already in classical
number theory as well as in the spectral theory of surfaces with cusps.
Our function $\zeta_{\mathbb{Z}}(s)$ can be seen to be an important
fudge factor (also called scattering determinant) for the Eisenstein
series, this comes from Selberg but he does not realize that it is
itself a spectral zeta function and writes just $\sqrt{\pi}\Gamma(s-1/2)/\Gamma(s)$.
We think that $\zeta_{\mathbb{Z}}(s)$ is undeniably present, the
question is whether its appearance is incidental or part of a general
structure. Evidence for the more structural picture could be that
for surfaces $\mathbb{Z}$ appears as the fundamental group of the
cusps. This leads to the question whether in higher dimensional Eisenstein
series, for example in Langlands' work \cite{La76}, other graph spectral
zeta functions occur. In this context, the \textbf{Langlands program,}
we observe that our function $\zeta_{\mathbb{Z}}(s)$ (or products
thereof) is essentially the value in the Bhanu-Murty-Gindikin-Karpelevich
formula at the archimedian place, see \cite[Ch. 8]{Bh60,FGKP16} related
to the Harish-Chandra $c$-function.
\begin{problem*}
Understand the role played by $\zeta_{\mathbb{Z}}(s)$, and perhaps
other spectral zeta functions, in the theory of Eisenstein series
and the Langlands program. For example, $\zeta_{\mathbb{Z}}(s)$ and
all the $p$-adic zeta functions in the last section might appear
together in the Langlands constant term, or Langlands $p$-adic Gindikin-Karpelevich
formula in \cite{La71}. 
\end{problem*}
The spectral zeta function of the finite cyclic graph $\mathbb{Z\mathrm{/}}n\mathbb{Z}$
(see e.g. \cite{CJK10,FK17} for details) is 
\[
\zeta_{\mathbb{Z\mathrm{/}}n\mathbb{Z}}(s)=\frac{1}{4^{s}}\sum_{k=1}^{n-1}\frac{1}{\sin^{2s}(\pi k/n)}.
\]

There exists some literature on finite trigonometric sums, which in
our context appears as the special values
\[
\zeta_{\mathbb{Z\mathrm{/}}n\mathbb{Z}}(m)
\]
for integral $m$. Some of these special values are recorded in \cite{FK17}.
In particular since we can transform $\sin^{-2}z$ to $\cot^{2}z$
we have special cases of \textbf{Dedekind sums} (which again leads
into the theory of Eisenstein series as I learnt from Claire Burrin),
see \cite{Z73}. In another context, in the first of the \textbf{Verlinde
formulas}, such sums appear, see \cite{Sz93,Z96} for mathematical
discussions. Let $\mathcal{N_{\mathrm{g,n,d}}}$ denote the moduli
space of semi-stable $n$-dimensional vector bundles over a fixed
Riemann surface of genus $g$ and having as determinant bundle a fixed
line bundle of degree $d$. The formula reads

\[
\dim_{\mathbb{C}}H^{0}(\mathcal{\mathcal{N_{\mathrm{g,2,0}}}\mathrm{,}L^{\mathrm{m}}\mathrm{)=\frac{(m+2)^{g-1}}{2^{g-1}}\sum_{k=1}^{m+1}\frac{1}{\sin^{2g-2}\frac{\pi k}{m+2}}.}}
\]
The right hand side can in our terminology be written as 
\[
(m+2)^{g-1}2^{g-1}\zeta_{\mathbb{Z\mathrm{/}}(m+2)\mathbb{Z}}(g-1).
\]
The lead term in the asymptotics as $m\rightarrow\infty$ was considered
by Witten \cite{W91} to evaluate the volume of the moduli space in
question and here the well-known special values $\zeta(2(g-1))$ of
the Riemann zeta function appears. This is consistent with the results
in \cite{FK17} and indeed earlier in a numerical analysis paper \cite{Si04}.
\begin{problem*}
Is the appearance of the special values of spectral zeta function
of cyclic graphs just a coincidence, or are there other cases of the
Verlinde formulas that allow interpretations as spectral zeta functions
of graphs? If so, this would be intriguing and demand for an explanation.
\end{problem*}
In \cite{FK17}, using some of the methods in \cite{CJK10,CJK12},
we study the asymptotics of the spectral zeta function of $\mathbb{Z\mathrm{^{d}/}}A_{n}\mathbb{Z}^{d}$
as $n\rightarrow\infty$ motivated in particular by statistical physics.
One sees there how in the asymptotics, the spectral zeta functions
of the infinite graphs and manifold spectral zeta functions appear.
In \cite{FK17} and \cite{F16} relations to analytic number theory
are discussed. Notably there are reformulations of the \textbf{Riemann
hypothesis }and the generalized Riemann hypothesis for certain Dirichlet
$L$-functions. For example in \cite{F16}: Let $m\geq3$ and let
$\chi$ be a primitive and even Dirichlet character modulo $m$. For
$n\geq1$, define the cyclic graph $L$-function 
\[
L_{n}(s,\chi):=\sum_{k=1}^{mn-1}\frac{\chi(k)}{\sin^{s}(\pi k/mn)}.
\]
Let $\Lambda_{n}(s,\chi)=n^{-s}(\pi/k)^{s/2}\Gamma(s/2)L_{n}(s,\chi).$
Recall the classical Dirichlet $L$-function
\[
L(s,\chi)=\sum_{k=1}^{\infty}\frac{\chi(k)}{k^{s}}.
\]
Friedli's theorem relates the Riemann hypothesis, on the location
of the zeros of this latter function to a functional relation for
the graph functions that imitates the well-known one for $L$ itself.
(This extension of \cite{FK17} seems to me rather surprising.) 
\begin{thm}
\cite{F16} The following two statements are equivalent:
\end{thm}

\begin{enumerate}
\item \emph{For all $s$ with $0<Re(s)<1$ and $Im(s)\geq8$ we have 
\[
\lim_{n\rightarrow\infty}\frac{\left|\Lambda_{n}(s,\chi)\right|}{\left|\Lambda_{n}(1-s,\overline{\chi})\right|}=1;
\]
}
\item \emph{In the region $0<Re(s)<1$ and $Im(s)\geq8$ , all zeros of
$L(s,\chi)$ have real part $1/2$. }
\end{enumerate}
The first statement holds in any case for all $s$ where $L$ does
not vanish. The appearance of a restriction on the imaginary part
has a substantial reason, see Lemma 3.1 in the proof of the theorem
in \cite{F16}. A related lemma with a similar restriction was proved
in the case of the Riemann zeta function and the restriction was shown
to be essential \cite{FK17}. 

\subsection{A few comments on further examples}

As it is pointed out in \cite{FK17} the spectral zeta functions of
such fundamental infinite graphs as the regular trees $T_{q+1}$ and
the standard lattices $\mathbb{Z}^{d}$ lead into \textbf{hypergeometric
functions} of several variables, more precisely, specializations of
these functions. 

The spectral zeta function of $\mathbb{Z}^{d}$ is
\[
\zeta_{\mathbb{Z}^{d}}(s)=\frac{d^{-s+1/2}}{\sqrt{2\pi}}\frac{\Gamma((s+1)/2)}{\Gamma(s)}F_{C}^{(d)}(s/2,(s+1)/2;1,1,...,1;1/d^{2},1/d^{2},...,1/d^{2}),
\]
where $F_{C}^{(d)}$ is one of the Lauricella hypergeometric functions
in $d$ variables.
\begin{problem*}
It is remarked in \cite[p. 49]{Ex76} that no integral representation
of Euler type has been found for $F_{C}$. We note that if one instead
of the heat kernel starts with the spectral measure in defining $\zeta_{\mathbb{Z}^{d}}(s)$,
we do get such an integral representation, at least for special parameters.
Does this lead to the missing Euler-type integral representation formula?
\end{problem*}
The $(q+1)$-regular tree (or Bethe lattice in physics parlance) is
the universal covering of $(q+1)$-regular graphs and is therefore
a fundamental graph. In \cite{FK17} an expression for the corresponding
zeta function is found, interestingly via an Euler-type integral that
Picard considered and which leads into Appell's hypergeometric function
$F_{1}$: 
\[
\zeta_{T_{q+1}}(s)=\frac{q(q+1)}{(q-1)^{2}(\sqrt{q}-1)^{2s}}F_{1}(3/2,s+1,1,3;u,v),
\]
with $u=-4\sqrt{q}/(\sqrt{q}-1)^{2}$ and $v=4\sqrt{q}/(\sqrt{q}+1)^{2}$.
\begin{problem*}
What functional equations do these spectral zeta functions have? In
view of the many symmetries that such hypergeometric functions have,
could one hope for an $s$ vs $1-s$ symmetry or similar identities?
\end{problem*}
The determinant of the Laplacian (of course removing the trivial eigenvalue
$0$ from the product) of a finite graph is known to count the number
of spanning trees (with a root) of the graph. This is called Kirchhoff's
matrix-tree theorem. For infinite graphs the corresponding determinant
is sometimes related to Mahler measures (in number theory) and more
generally to \textbf{Fuglede-Kadison determinants} (in operator algebras).
See \cite{Ly10,CJK10,CJK12} for more about these connections. A famous
value is the determinant of the graph $\mathbb{Z}^{2}$ which is $4/\pi$
times Catalan's constant; one my wonder if there are other such values
in terms of special values of Dirichlet $L$-functions.

\section{Totally disconnected case}

\subsection{$\mathbb{Q}_{p}$ and $\mathbb{Q}_{p}/\mathbb{Z}_{p}$ }

This follows \cite{CZ17} and the recent master thesis of Mårten Nilsson
\cite{Ni18}, see these two sources for further references. Let $\mathbb{Q}_{p}$
denote the $p$-adic numbers and $\mathbb{Z}_{p}$ the $p$-adic integers.
Let $dy$ be the Haar measure on the locally compact additive group
$\mathbb{Q}_{p}$ normalized so that the measure of $\mathbb{Z}_{p}$
is $1$. The Taibleson-Vladimirov Laplcaian can be defined like pseudo-differential
operators via the Fourier transform, alternatively it is explicitly
given as an integral as follows:
\[
\Delta f(x)=\frac{p^{2}-1}{1-p^{-3}}\int_{\mathbb{Q}_{p}}\frac{f(x)-f(y)}{\left|x-y\right|_{p}^{3}}dy
\]
for suitable functions $f:\mathbb{Q}_{p}\rightarrow\mathbb{C}$. This
gives rise in the usual ways to a heat equation and a heat kernel,
which in this case turns out to be:
\[
K_{p}(x,t)=\sum_{k=-\infty}^{\infty}\left(e^{-tp^{2k}}-e^{-tp^{2k+2}}\right)p^{k}C_{p^{-k}}(x),
\]
where $C_{p^{n}}$ denotes the characteristic function of the the
ball $B_{p^{n}}=\left\{ x\in\mathbb{Q}_{p}:\left|x\right|_{p}\leq p^{n}\right\} $.
There is another formula for this function:
\[
K_{p}(x,t)=\sum_{m=0}^{\infty}\frac{(-1)^{m}t^{m}}{m!}\frac{1-p^{2m}}{1-p^{-2m-1}}\left|x\right|_{p}^{-2m-1},
\]
valid for $x\neq0$. Now passing to the quotient $\mathbb{Q}_{p}/\mathbb{Z}_{p}$
(which is analogous to the periodization done above and is here an
integral over $\mathbb{Z}_{p}$), and after that taking the Mellin
transform dividing by $\Gamma(s)$, leads to the corresponding spectral
zeta function which is 
\[
\zeta_{p}(s)=(1-p^{-1})\frac{p^{1-2s}}{1-p^{1-2s}}=\frac{p^{1-2s}-p^{-2s}}{1-p^{1-2s}}=\frac{p-1}{p^{2s}-p}.
\]
Here we complete this function in the following manner:
\[
\xi_{p}(s)=\sin(2\pi s)p^{s}\zeta_{p}(s),
\]
then we obtain a functional equation of the usual type in the most
symmetric form:

\[
\xi_{p}(1-s)=\xi_{p}(s).
\]
Since the zeta functions became so simple this relation is more trivial
than what we saw in the other contexts. But still, we see a pattern:
From the line $\mathbb{R}$ and the circle $\mathbb{R}/\mathbb{Z}$,
to the graphs $\mathbb{Z}$ and $\mathbb{Z\mathrm{/n\mathbb{Z}}}$,
and now $\mathbb{Q}_{p}$ and $\mathbb{Q}_{p}/\mathbb{Z}_{p}$, their
associated spectral zeta functions have a non-obvious symmetry $s$
vs $1-s$. 

In another direction taking the Laplace transform of the two heat
kernel expressions, Nilsson derives an identity valid for $Re(s)>(p/\left|x\right|_{p})^{2}$:
\[
\frac{1}{s\left|x\right|_{p}}\sum_{m=0}^{\infty}\frac{(-1)^{m}}{s^{m}\left|x\right|_{p}^{2m}}\frac{1-p^{2m}}{1-p^{-2m-1}}=\sum_{k=-\infty}^{\infty}\frac{p^{2}-1}{(1+p^{2k}s)(p^{2}+p^{2k}s)}p^{k}C_{p^{k}}(x).
\]
Specializing to certain $x$ gives an identity without $p$-adics,
but only involving ordinary integers. For example, with $x=1$ and
$s>p^{2}$,
\[
\frac{1}{s}\sum_{m=0}^{\infty}\frac{(-1)^{m}}{s^{m}}\frac{1-p^{2m}}{1-p^{-2m-1}}=\sum_{k=0}^{\infty}\frac{p^{2}-1}{(1+p^{2k}s)(p^{2}+p^{2k}s)}\cdot p^{k}.
\]

It is natural to compare these considerations with the celebrated
\textbf{Tate's thesis}, where Tate in particular after a Mellin transform,
obtains for each $p$ the local Euler factor $(1-p^{-s})^{-1}$, instead
of our $\zeta_{p}(s)$. This brings us back to the first paragraph
of this paper.

Section de mathématiques, Université de Genève, 2-4 Rue du Lièvre,
Case Postale 64, 1211 Genève 4, Suisse 

e-mail: anders.karlsson@unige.ch 

and

Matematiska institutionen, Uppsala universitet, Box 256, 751 05 Uppsala,
Sweden 

e-mail: anders.karlsson@math.uu.se
\end{document}